# A REMARK ON THE GENUS OF THE INFINITE QUATERNIONIC PROJECTIVE SPACE

DONALD YAU

ABSTRACT. It is shown that only countably many spaces in the genus of $\mathbf{HP}^\infty$, the infinite quaternionic projective space, can admit essential maps from $\mathbf{CP}^\infty$, the infinite complex projective space. Examples of countably many homotopically distinct spaces in the genus of $\mathbf{HP}^\infty$ which admit essential maps from $\mathbf{CP}^\infty$ are constructed. These results strengthen a theorem of McGibbon and Rector which states that among the uncountably many homotopy types in its genus, $\mathbf{HP}^\infty$ is the only one which admits a maximal torus.

## 1. INTRODUCTION

In an attempt to understand Lie group theory through homotopy theory, Rector [8] classified all the spaces in the genus of $\mathbf{HP}^\infty$, the infinite projective space over the quaternions. A space $X$ (or its homotopy type) is in the genus of $\mathbf{HP}^\infty$, denoted $X \in \mathrm{Genus}(\mathbf{HP}^\infty)$, if the $p$-localizations of $X$ and $\mathbf{HP}^\infty$ are homotopy equivalent for each prime $p$. Rector proved the following classification result in [8]:

**Theorem 1.1** (Rector). *Let $X$ be a space in the genus of $\mathbf{HP}^\infty$. Then for each prime $p$ there exists an invariant $(X/p) \in \{\pm 1\}$ so that the following holds:*

1. *The $(X/p)$ for $p$ primes provide a complete list of classification invariants for the genus of $\mathbf{HP}^\infty$.*
2. *Any combination of values of the $(X/p)$ can occur. In particular, the genus of $\mathbf{HP}^\infty$ is uncountable.*
3. *$(\mathbf{HP}^\infty/p) = 1$ for all primes $p$.*
4. *$X$ has a maximal torus if and only if $X$ is homotopy equivalent to $\mathbf{HP}^\infty$.*

Actually, for the last statement about the maximal torus, Rector proved it for the odd primes (i.e. that $X$ has a maximal torus implies $(X/p) = 1$ for all odd primes $p$), and then McGibbon [6] proved it for the prime $p = 2$ as well. Here $X$ is said to have a maximal torus if there exists a map from $BS^1 = \mathbf{CP}^\infty$ to $X$ whose homotopy theoretic fiber has the homotopy type of a finite complex.

For a space $X$ in the genus of $\mathbf{HP}^\infty$ which is not homotopy equivalent to $\mathbf{HP}^\infty$, the nonexistence of a maximal torus does not rule out the possibility that there could be some essential (i.e. non-nullhomotopic) maps from $\mathbf{CP}^\infty$ to $X$, and this can actually happen (see Example 1.3). The main purpose of this note is to show that this can (and it does) happen to only countably many spaces in the genus of $\mathbf{HP}^\infty$:

**Theorem 1.2.** *There are at most countably many spaces in the genus of $\mathbf{HP}^\infty$ which admit essential maps from $\mathbf{CP}^\infty$.*

The following example provides countably many homotopically distinct spaces in the genus of $\mathbf{HP}^\infty$ which admit essential maps from $\mathbf{CP}^\infty$.

**Example 1.3.** Let $p$ be any odd prime. Denote by $X(p)$ the space in the genus of $\mathbf{HP}^\infty$ with Rector invariants
$$(X(p)/q) = \begin{cases} 1 & \text{if } q \neq p, \\ -1 & \text{if } q = p. \end{cases}$$
Then there exists an essential map from $\mathbf{CP}^\infty$ to $X(p)$.

A couple of remarks are in order.

**Remark 1.4.** Our main result above (Theorem 1.2) can be regarded as an attempt to understand the set of homotopy classes of maps from $X$ to $Y$, where $X \in \text{Genus}(BG)$ and $Y \in \text{Genus}(BK)$ with $G$ and $K$ some connected compact Lie groups. This problem was posed and studied by Ishiguro, Møller, and Notbohm [5].

**Remark 1.5.** Theorem 1.2 is closely related to the nonexistence of embeddings of integral $K$-theory of finite $H$-spaces into $K^*(BT^n)$, analogous to the Adams-Wilkerson embedding of unstable algebras over the Steenrod algebra into the mod $p$ cohomology of $BT^n$. Adams once constructed a space $Z \in \text{Genus}(\mathbf{HP}^\infty)$ with the property that every nonzero element in $\overline{K}^0(Z)$ has infinitely many nonzero Chern classes; see [1, p. 79]. Therefore, $K(Z)$ cannot be embedded into $K(\mathbf{CP}^\infty)$ as a sub-$\lambda$-ring. Together with the result of Notbohm and Smith [7, Thm. 5.2] (or Dehon and Lannes [4]) our Theorem 1.2 shows that, in fact, there are only countably many spaces $X$ in the uncountable set Genus($\mathbf{HP}^\infty$) for which $K(X)$ can be embedded into $K(\mathbf{CP}^\infty)$ as sub-$\lambda$-rings.

The proof of Theorem 1.2 makes use of the results of Notbohm and Smith mentioned above. In fact, it follows from their result [7, Thm. 5.2] that two maps $f, g\colon \mathbf{CP}^\infty \to X \in \text{Genus}(\mathbf{HP}^\infty)$ are homotopic if and only if they have the same degree, where the degree $\deg(f)$ of such a map $f$ is defined by the equation

(1.6) $\qquad K^*(f)(b^2 u_X) = \deg(f)(b^2\lambda^2) + \text{higher terms in } b\lambda$

in complex $K$-theory. Here $b^2 u_X \in K(X)$ and $b\lambda \in K(\mathbf{CP}^\infty)$ are ring generators, and $b \in K^{-2}(\text{pt})$ is the Bott element. Therefore, to prove Theorem 1.2 it suffices to establish the following lemma.

**Lemma 1.7.** *Let $p$ be an odd prime and $k$ be a non-zero integer relatively prime to $p$. Suppose that for some space $X$ in the genus of $\mathbf{HP}^\infty$ there exists a map $f\colon \mathbf{CP}^\infty \to X$ of degree $k$. Then $(k/p) = (X/p)$.*

Here $(X/p)$ is the Rector invariant of $X$ as defined in [8], and $(k/p)$ is the Legendre symbol of $k$ (i.e. it is either $1$ or $-1$ according as whether $k$ is a quadratic residue modulo $p$ or not). It is an elementary result in number theory that if $k$ is relatively prime to $p$, then $(k/p) \equiv k^{(p-1)/2}$ mod $p$.

It follows from Lemma 1.7 that for each non-zero integer $k$, there are at most finitely many spaces in the genus of $\mathbf{HP}^\infty$ which admit a map of degree $k$ from $\mathbf{CP}^\infty$. Theorem 1.2 then follows from this and the result of Notbohm and Smith mentioned above.

The proof of Lemma 1.7 is made possible by the result in [10] which identifies the filtered ring $K(X)$ with $K(\mathbf{HP}^\infty)$ (up to a non-canonical filtered ring isomorphism) for any $X$ in the genus of $\mathbf{HP}^\infty$, and which also computes the Adams operations of $K(X)$ modulo certain ideals (see Proposition 2.1).

**Remark 1.8.** For certain integers $k$, one can be more specific about how many spaces there could be in the genus of $\mathbf{HP}^\infty$ which admit maps from $\mathbf{CP}^\infty$ of degree $k$. For example, with some calculations similar to those in §3 involving the Adams operations $\psi^2$, one can show that there does not exist any map from $\mathbf{CP}^\infty$ to $X$ of degree $2$ for any space $X$ in the genus of $\mathbf{HP}^\infty$. However, such calculations have not led to any substantial strengthening of Theorem 1.2, and therefore we omit them.

The rest of this note is organized as follows. In §2 we recall some results from [10] (with sketches of proofs) which are needed in §3, in which the proof of Lemma 1.7 is given. Example 1.3 is proved in §4. For use in §2, this note ends with an appendix about localization towers.

## 2. $K$-theory filtered ring

In this section, special cases of some results from [10] concerning the $K$-theory filtered rings of spaces in the genus of $\mathbf{HP}^\infty$ are reviewed.

To begin with, a filtered ring is a pair $(R, \{I_n^R\})$ consisting of: (1) a commutative ring $R$ with unit, and (2) a decreasing filtration $R = I_0^R \supset I_1^R \supset \cdots$ of ideals of $R$. A map between two filtered rings is a ring homomorphism which preserves the filtrations. With these maps as morphisms, the filtered rings form a category.

Every space $Z$ of the homotopy type of a CW complex gives rise naturally to an object $(K(Z), \{K_n(Z)\})$ in this category. Here $K(Z)$ and $K_n(Z)$ denote, respectively, the complex $K$-theory ring of $Z$ and the kernel of the restriction map $K(Z) \to K(Z_{n-1})$, where $Z_{n-1}$ denotes the $(n-1)$-skeleton of $Z$. Using a different CW structure of $Z$ will not change the filtered ring isomorphism type of $K(Z)$, as can be easily seen by using the cellular approximation theorem. The symbol $K_s^r(Z)$ denotes the quotient $K^r(Z)/\ker(K^r(Z) \to K^r(Z_{s-1}))$.

The following result, which is proved in [10], will be needed in §3. For the reader's convenience we include here a sketch of the proof. In what follows $b \in K^{-2}(\mathrm{pt})$ will denote the Bott element.

**Proposition 2.1.** *Let $X$ be a space in the genus of $\mathbf{HP}^\infty$. Then the following statements hold:*

1. *There exists an element $u_X \in K_4^4(X)$ such that $K(X) = \mathbf{Z}[[b^2 u_X]]$ as a filtered ring.*
2. *For any odd prime $p$, the Adams operation $\psi^p$ satisfies*

$$(2.2) \quad \psi^p\left(b^2 u_X\right) = \left(b^2 u_X\right)^p + 2\left(X/p\right) p \left(b^2 u_X\right)^{(p+1)/2} + pw + p^2 z,$$

   *where $w$ and $z$ are some elements in $K_{2p+3}^0(X)$ and $K_4^0(X)$, respectively. In particular, we have*

$$(2.3) \quad \psi^p\left(b^2 u_X\right) = 2\left(X/p\right) p \left(b^2 u_X\right)^{(p+1)/2} \quad (mod\ K_{2p+3}^0(X)\ and\ p^2).$$

*Sketch of the proof of Proposition 2.1.* For the first assertion, Wilkerson's proof of the classification theorem [9, Thm. I] of spaces of the same $n$-type for all $n$ can be easily adapted to show the following (see Proposition A.1 in Appendix A): There is a bijection between the pointed set of isomorphism classes of filtered rings $(R, \{I_n^R\})$ with the properties, (a) $R \xrightarrow{\cong} \varprojlim_n R/I_n^R$ and (b) $R/I_n^R \cong K(\mathbf{HP}^\infty)/K_n(\mathbf{HP}^\infty)$ as filtered rings, and the pointed set $\varprojlim_n^1 \mathrm{Aut}(K(\mathbf{HP}^\infty)/K_n(\mathbf{HP}^\infty))$. Here $\mathrm{Aut}(-)$ denotes the group of filtered ring automorphisms, and the $\varprojlim^1$ of a tower of not-necessarily abelian groups is as defined in [3]. It is not difficult to check that the hypothesis on $X$ implies that the filtered ring $K(X)$ has the above two properties, (a) and (b). Moreover, by analyzing the subquotients $K_n(\mathbf{HP}^\infty)/K_{n+1}(\mathbf{HP}^\infty)$, one can show that the map

$$\mathrm{Aut}(K(\mathbf{HP}^\infty)/K_{n+1}(\mathbf{HP}^\infty)) \to \mathrm{Aut}(K(\mathbf{HP}^\infty)/K_n(\mathbf{HP}^\infty))$$

is surjective for each $n > 4$. The point is that any automorphism of $K(\mathbf{HP}^\infty)/K_n(\mathbf{HP}^\infty)$ can be lifted to an endomorphism of $K(\mathbf{HP}^\infty)/K_{n+1}(\mathbf{HP}^\infty)$ without any difficulty. Then, since the quotient $K(\mathbf{HP}^\infty)/K_{n+1}(\mathbf{HP}^\infty)$ is a finitely generated abelian group, one only has to observe that the chosen lift is surjective. Now it follows that the above $\varprojlim^1$ is the one-point set, and hence $K(X) \cong K(\mathbf{HP}^\infty)$ as filtered rings. This establishes the first assertion.

The second assertion concerning the Adams operations $\psi^p$ is an immediate consequence of Atiyah's theorem [2, Prop. 5.6] and the definition of the invariants $(X/p)$ given by Rector [8]. □

## 3. Proof of Lemma 1.7

In this section the proof of Lemma 1.7 is given.

First, recall that the complex $K$-theory of $\mathbf{CP}^\infty$ as a filtered $\lambda$-ring is given by $K(\mathbf{CP}^\infty) = \mathbf{Z}[[b\lambda]]$ for some $\lambda \in K_2^2(\mathbf{CP}^\infty)$. The Adams operations on the generator are given by

(3.1) $\qquad \psi^r(b\lambda) = (1+b\lambda)^r - 1 \quad (r = 1, 2, \dots).$

From now on $k$, $p$, $f$, and $X$ are as in Lemma 1.7, i.e. $k$ is a fixed non-zero integer, $p$ is an odd prime which is relatively prime to $k$, $X$ is some space in the genus of $\mathbf{HP}^\infty$, and $f$ is some (essential) map $\mathbf{CP}^\infty \to X$ of degree $k$. Denote by $f^*$ the map induced by $f$ in $K$-theory.

We will compare the coefficients of $(b\lambda)^{p+1}$ in the equation

(3.2) $\quad f^*\psi^p\left(b^2 u_X\right) = \psi^p f^*\left(b^2 u_X\right) \quad (\mathrm{mod}\ K_{2p+3}^0(\mathbf{CP}^\infty) \text{ and } p^2).$

For the left-hand side of Eq. (3.2), it follows from Eq. (1.6) and (2.3) that, modulo $K_{2p+3}^0(\mathbf{CP}^\infty)$ and $p^2$,

(3.3)
$$f^*\psi^p\left(b^2 u_X\right) = 2\,(X/p)\,p\,\left(kb^2\lambda^2\right)^{(p+1)/2} = 2\,(X/p)\,p\,k^{(p+1)/2}(b\lambda)^{p+1}.$$

Similarly, for the right-hand side of Eq. (3.2), it follows from Eq. (1.6) and (3.1) (when $r = 2$) that, modulo $K^0_{2p+3}(\mathbf{CP}^\infty)$ and $p^2$,

(3.4) $\quad \psi^p f^* \left( b^2 u_X \right) = k\psi^p(b^2\lambda^2) = k\psi^p(b\lambda)^2 = 2pk(b\lambda)^{p+1}.$

Combining Eq. (3.3) and (3.4) we obtain the congruence relation

$$2\left(X/p\right) p\, k^{(p+1)/2} \equiv 2pk \pmod{p^2},$$

which is equivalent to

(3.5) $\quad (X/p)\,(k/p) \equiv (X/p)\, k^{(p-1)/2} \equiv 1 \pmod{p}$

because $p$ is odd and is relatively prime to $k$. Hence $(X/p) = (k/p)$, as desired.

This finishes the proof of Lemma 1.7.

## 4. Proof of Example 1.3

Fix an odd prime $p$ and let $X(p) \in \text{Genus}(\mathbf{HP}^\infty)$ be as in Example 1.3. To construct an essential map $f\colon \mathbf{CP}^\infty \to X(p)$, first recall that $X(p)$ is obtained as the homotopy inverse limit [3]

$$X(p) = \text{holim}_q \left\{ \mathbf{HP}^\infty_{(q)} \xrightarrow{r_q} \mathbf{HP}^\infty_{(0)} \xrightarrow{n_q} \mathbf{HP}^\infty_{(0)} \right\}.$$

Here $\mathbf{HP}^\infty_{(q)}$ (resp. $\mathbf{HP}^\infty_{(0)}$) denotes the $q$-localization (resp. rationalization) of $\mathbf{HP}^\infty$, $q$ runs through all primes, $r_q$ is the rationalization map. The integers $n_q$ are chosen as follows:

$$n_q = \begin{cases} 1 & \text{if } q \neq p \\ k & \text{if } q = p, \end{cases}$$

where $k$ is an integer with $1 < k < p$ and is a quadratic non-residue modulo $p$.

Denote by $l_q\colon \mathbf{HP}^\infty \to \mathbf{HP}^\infty_{(q)}$ the $q$-localization map and by $i\colon \mathbf{CP}^\infty \to \mathbf{HP}^\infty$ the map induced by the standard inclusion $S^1 \to S^3$. (We are identifying $BS^1$ with $\mathbf{CP}^\infty$ and $BS^3$ with $\mathbf{HP}^\infty$.) Define maps $f_q\colon \mathbf{CP}^\infty \to \mathbf{HP}^\infty_{(q)}$ for $q$ primes by

$$f_q = \begin{cases} l_q \circ i \circ Bk & \text{if } q \neq p, \\ l_p \circ i & \text{if } q = p. \end{cases}$$

Here $Bk\colon \mathbf{CP}^\infty \to \mathbf{CP}^\infty$ is the map induced by multiplication by $k$ in $S^1$. It is then easy to see that $n_q \circ r_q \circ f_q = n_{q'} \circ r_{q'} \circ f_{q'}$ for any two primes $q$ and $q'$. Therefore, the maps $f_q$ glue together to yield a map

$$f\colon \mathbf{CP}^\infty \to X(p)$$

through which every map $f_q$ factors. This map $f$ is easily seen to be essential, since the composite

$$k \circ r_p \circ f_p \colon \mathbf{CP}^\infty \to \mathbf{HP}^\infty_{(0)} = K(\mathbf{Q}, 4)$$

represents $k$ times the square of the rational cohomology generator of $\mathbf{CP}^\infty$ and hence is essential.

This proves Example 1.3.

## Appendix A. Localization towers

In this appendix we describe a result about localization towers which was used in §2.

Let $\mathcal{M}$ be a pointed model category. A *localization functor* on $\mathcal{M}$ is a pair $(L, \eta)$ consisting of: (1) an endofunctor $L \colon \mathcal{M} \to \mathcal{M}$, and (2) a natural transformation $\eta \colon \mathrm{Id} \to L$ from the identity functor to $L$, such that for every object $X$ in $\mathcal{M}$, the two maps $L\eta_X, \eta_{LX} \colon LX \to L^2 X$ are homotopic and are both weak equivalences. The reference to $\eta$ is often omitted and one speaks of $L$ as a localization functor. In the literature, a localization functor is sometimes referred to as a *homotopy idempotent monad*.

A *localization tower* on $\mathcal{M}$ consists of: (1) localization functors $L_n$ ($n = 0, 1, 2, \ldots$) on $\mathcal{M}$, and (2) natural transformations $L_{n+1} \to L_n$ induced by the universal property of $L_{n+1}$-localization. For example, the Postnikov tower $\{P_{S^n}\}$ and the stable chromatic tower $\{L_n\} = \{L_{E(n)}\}$ are both localization towers.

A special case of the following classification result was used in §2.

**Proposition A.1.** *Let $\mathcal{M}$ be a pointed model category and $\{L_n\}$ be a localization tower on $\mathcal{M}$. Let $X$ be an object in $\mathcal{M}$ such that the natural map $X \to \mathrm{holim}_n L_n X$ is a weak equivalence. Then there exists a bijection between the following two sets:*

1. *The set of all weak equivalence classes $[Y]$ of objects in $\mathcal{M}$ for which $L_n Y$ and $L_n X$ are weakly equivalent for all $n \geq 0$, and the natural map $Y \to \mathrm{holim}_n L_n Y$ is a weak equivalence.*
2. $\varprojlim{}^1 \mathrm{Aut}(L_n X)$.

Here $\mathrm{Aut}(Z)$ denotes the group of homotopy classes of self-weak-equivalences of the functorial fibrant-cofibrant replacement of $Z$.

This proposition is a straightforward generalization of Wilkerson's proof of his classification theorem of spaces of the same $n$-type for all $n$ [9].

In §2 (and also [10]) this result is applied to the category of (commutative) filtered rings with the trivial model category structure (i.e.

weak equivalences are just the isomorphisms and every map is both a fibration and a cofibration), $L_n((R, \{I_n^R\})) = R/I_n^R$ with the induced quotient filtration, and $X = (K(\mathbf{HP}^\infty), \{K_n(\mathbf{HP}^\infty)\})$.

## Acknowledgment

Part of the materials in this short paper was presented in a talk in the 2001 International Conference in Algebraic Topology (Isle of Skye, Scotland). The author would like to thank the organizers of that conference for their wonderful arrangements. Thanks also go to Dietrich Notbohm for his interest and comments in this work. Finally, the author expresses his sincerest gratitude to his advisor Haynes Miller for many hours of inspiring conversations and for his bottomless patience with the author's mistakes.

Department of Mathematics, MIT, Cambridge, Massachusetts, USA
*E-mail address*: donald@math.mit.edu